# The size Ramsey number of a directed path

Ido Ben-Eliezer[*]    Michael Krivelevich [†]    Benny Sudakov [‡]

May 25, 2010

**Abstract**

Given a graph $H$, the size Ramsey number $r_e(H, q)$ is the minimal number $m$ for which there is a graph $G$ with $m$ edges such that every $q$-coloring of $G$ contains a monochromatic copy of $H$. We study the size Ramsey number of the directed path of length $n$ in oriented graphs, where no antiparallel edges are allowed. We give nearly tight bounds for every fixed number of colors, showing that for every $q \geq 1$ there are constants $c_1 = c_1(q), c_2$ such that

$$\frac{c_1(q) n^{2q} (\log n)^{1/q}}{(\log \log n)^{(q+2)/q}} \leq r_e(\overrightarrow{P_n}, q+1) \leq c_2 n^{2q} (\log n)^2.$$

Our results show that the path size Ramsey number in oriented graphs is asymptotically larger than the path size Ramsey number in general directed graphs. Moreover, the size Ramsey number of a directed path is polynomially dependent in the number of colors, as opposed to the undirected case.

Our approach also gives tight bounds on $r_e(\overrightarrow{P_n}, q)$ for general directed graphs with $q \geq 3$, extending previous results.

## 1 Introduction

Given an integer $q > 0$, we write $G \to (H)_q$ if every $q$-coloring of $E(G)$ contains a monochromatic copy of $H$.

The study of size Ramsey numbers (initiated in [8]) is concerned with the following questions. Given a graph $H$, what is the minimum number of edges $m$ for which there is a graph $G$ with $m$ edges such that $G \to (H)_q$?

Denote by $r_e(H, q)$ the size Ramsey number of $H$ with respect to coloring with $q$ colors. That is,

$$r_e(H, q) = \min\{|E(G)| \mid G \to (H)_q\}.$$

The study of $r_e(K_n, q)$ is essentially equivalent to the study of the original Ramsey number. Namely, it can be verified that if $r_e(K_n, q) = m$ then a clique with exactly $m$ edges has the desired property. This result is attributed to Chvátal in [8].

[*]School of Computer Science, Raymond and Beverly Sackler Faculy of Exact Sciences, Tel Aviv University, Tel Aviv 69978, Israel, e-mail: **idobene@post.tau.ac.il**. Research supported in part by an ERC advanced grant.

[†]School of Mathematical Sciences, Raymond and Beverly Sackler Faculty of Exact Sciences, Tel Aviv University, Tel Aviv 69978, Israel, e-mail:   **krivelev@post.tau.ac.il**. Research supported in part by USA-Israel BSF grant 2006322, by grant 1063/08 from the Israel Science Foundation, and by a Pazy memorial award.

[‡]Department of Mathematics, UCLA, Los Angeles, CA 90095. Email: **bsudakov@math.ucla.edu**. Research supported in part by NSF CAREER award DMS-0812005 and by a USA-Israeli BSF grant.

Beck proved [6] that for every constant $q > 0$, the size Ramsey number of the undirected path on $n$ vertices is linear in $n$, answering a question of Erdős. Namely, there is a constant $c_q$ such that

$$r_e(P_n, q) \leq c_q n.$$

Later Alon and Chung [3] provided an explicit construction of graphs with this property. The linearity of the size Ramsey number for bounded degree trees was proved by Friedman and Pippenger [10].

**Directed graphs.** In this work we focus on the size Ramsey number in directed graphs, and in particular on the size Ramsey number of a directed path on $n$ vertices.

Raynaud [14] proved that every red-blue coloring of the complete symmetric directed graph (a directed graph in which between every two vertices there are edges in both directions) has a Hamilton cycle which is the union of two monochromatic paths (a simple proof can be found in [12]). In particular, this shows that the size Ramsey number of the path $\overrightarrow{P_n}$ for directed graphs with antiparallel edges is $O(n^2)$. Reimer [15] proved that every digraph with the property that every red-blue coloring has a path of length $n$ must have $\Omega(n^2)$ edges, therefore proving that Raynaud's bound for non-simple directed graphs is tight up to a constant factor.

In this work we provide nearly tight bounds for the size Ramsey number of directed paths in oriented graphs, where no parallel edges are allowed. We show that this number is asymptotically larger than the path size Ramsey number when such edges are allowed. Our approach also generalizes to the case of non-simple directed graphs, and we provide nearly tight bounds for such graphs for every constant $q \geq 3$.

Our first result is a lower bound on the size Ramsey number.

**Theorem 1.** *For every $q \geq 1$ there is a constant $c_1 = c_1(q)$ such that*

$$r_e(\overrightarrow{P_n}, q+1) \geq \frac{c_1 n^{2q} (\log n)^{1/q}}{(\log \log n)^{(q+2)/q}}.$$

We also have the following almost matching upper bound.

**Theorem 2.** *There is an absolute constant $c_2$ such that for every $q \geq 1$*

$$r_e(\overrightarrow{P_n}, q+1) \leq c_2 n^{2q} \cdot (\log n)^2.$$

*Moreover, a random tournament $T_N$ on $N = \Theta(n^q \log n)$ vertices satisfies $T_N \to (\overrightarrow{P_n})_{q+1}$ with high probability.*

We stress that Theorem 1 proves that the size Ramsey number of a directed path in oriented graphs is asymptotically larger than the size Ramsey number in general directed graphs for any fixed number of colors.

In the course of the proof of Theorem 2, we actually prove the following asymmetric Ramsey property.

**Theorem 3.** *There is an absolute constant $c$ such that the following holds. There is an oriented graph $G$ with $cn^2 \cdot (\log n)^2$ edges such that every red-blue coloring of $E(G)$ contains a red path of length $n$ or a blue path of length $n \log n$.*

Our approach also gives matching lower and upper bounds for the case of directed non-simple graphs with more than 2 colors, where antiparallel edges are allowed. This generalizes the results of Raynaud [14] and Reimer [15] to any fixed $q$.

**Proposition 4.** *The size Ramsey number of a directed path on $n$ vertices for $q+1$ colors in directed, non-simple graphs is $\Omega(n^{2q})$.*

**Proposition 5.** *The size Ramsey number of a directed path on $n$ vertices for $q+1$ colors in directed, non-simple graphs is $O(n^{2q})$.*

Note that here we show that the size Ramsey number (both for simple or non-simple directed graphs) of a directed path is polynomially dependent in the number of colors, in contrast to the undirected case where changing the number of colors changes the size Ramsey number by a constant factor, as follows from the above mentioned result of Beck [6].

**Organization of the paper.** The rest of the paper is organized as follows. In Section 2 we provide some notation and definitions and present well known theorems that will be used later on. In Section 3 we present the proofs of our lower bounds, proving Theorem 1 and Proposition 4. In Section 3 we prove the upper bounds, showing Theorem 2, Theorem 3 and Proposition 5. Finally, in Section 5 we present some concluding remarks and open problems.

Throughout the paper we assume that the underlying parameter $n$ is large enough. We do not try to optimize constants, and all logarithms are in base 2. We ignore all floor and ceiling signs whenever these are not crucial.

## 2 Preliminaries

All the graphs we consider are directed. An oriented (or simple) graph is a graph with no antiparallel (or opposite) edges. Given a set of vertices $A \subseteq V$, we denote by $E(A)$ the set of edges in the induced subgraph $G[A]$. The in-degree of a vertex $v$, denoted by $d^-(v)$, is the number of edges that are directed into $v$, and the out-degree of $v$, denoted by $d^+(v)$ is the number of edges directed from $v$. The degree of $v$ is $d^+(v) + d^-(v)$, and we let $\Delta(G)$ be the maximum degree in a graph $G$ and $\delta(G)$ be the minimum degree in $G$. For a vertex $v$ we let $N^+(v) = \{u \in V : (v, u) \in E\}$ and call this set the *out-neighbors* of $v$. We also let $N^-(v) = \{u \in V : (u, v) \in E\}$ and call this set the *in-neighbors* of $v$. For a set of vertices $A$ we let $N^+(A)$ be $\bigcup_{a \in A} N^+(a)$. A set of vertices is *acyclic* if it does not span a cycle. The *length* of a directed path is the number of edges it contains. The *edge density* of a directed graph $G = (V, E)$ is $\frac{|E|}{|V|^2}$.

The complete graph on $n$ vertices, denoted by $K_n$, is an undirected graph for which every pair of vertices are connected. A *complete symmetric* directed graph is a non-simple directed graph where between every two vertices there are edges in both directions. A tournament is an oriented graph where between every two vertices there is an edge in exactly one of the directions.

A $k$-coloring of a graph is a mapping of the vertices to $\{1, \ldots, k\}$ such that every two adjacent vertices are mapped to distinct values. The chromatic number of a graph $G$, denoted by $\chi(G)$, is the minimum $k$ such that $G$ is $k$-colorable. It is well known that every graph of maximum degree $d$ is $(d+1)$-colorable. A Hamilton cycle in $G$ is a cycle that visits every vertex in $G$ exactly once.

We will use the following two theorems several times throughout the paper. The first one, also mentioned in the introduction, is due to Raynaud [14] (see, e.g., [12] for a proof).

**Theorem 2.1.** *Let $G$ be a complete symmetric directed graph on $t$ vertices. Then every red-blue coloring of $E(G)$ contains a Hamilton cycle which is the union of two monochromatic paths. In particular, it contains a monochromatic path of length $t/2$.*

We also need the following simple theorem that was proved independently by Gallai [11] and Roy [16] (see, e.g., [13], Chapter 9, Problem 9, for a proof).

**Theorem 2.2.** *Let $G$ be a directed graph with no path longer than $t$. Then $G$ is $(t+1)$-colorable.*

## 3 Lower bounds

In this section we show that every oriented graph with relatively few edges has an edge coloring without a long monochromatic path. We first prove that every sparse oriented graph has a large acyclic set, and then use it to show that every graph admits a partition into a relatively small number of independent sets and acyclic sets. We conclude by showing that given such a partition we can color the edges with no long monochromatic path.

### 3.1 Sparse graphs have large acyclic sets

The following lemma is well known and easy.

**Lemma 3.1.** *Every tournament on $n$ vertices has an acyclic set of size $\lfloor \log n \rfloor + 1$.*

*Proof.* We may assume that $n = 2^k$ for some integer $k$, as otherwise we can take any subotournament of size $2^{\lfloor \log n \rfloor}$. We prove it by induction on $k$, and note that the case $k = 0$ is trivial. Suppose that the induction hypothesis is true for $k$, and we next prove it for $k+1$. Indeed, every tournament on $2^{k+1}$ vertices contains a vertex $v$ with out-degree at least $2^k$. By induction hypothesis, $N^-(v)$ contains an acyclic set $A$ of size $k+1$, and thus $A \cup \{v\}$ is an acyclic set of size $k+2$, as required. □

Since every oriented graph is a subgraph of a tournament, we get the following direct consequence.

**Corollary 3.2.** *Every oriented graph on $n$ vertices has an acyclic set of size $\log n$.*

Our main result in this subsection is that every sparse oriented graph has a larger acyclic set. This is a directed version of a lemma by Erdős and Szemerédi [9].

**Lemma 3.3.** *There is an absolute constant $c > 0$ such that the following holds. Every oriented graph $G$ with $n$ vertices and at most $\varepsilon n^2$ edges contains an acyclic set of size $\frac{c \log n}{\varepsilon \log (1/\varepsilon)}$.*

*Proof.* Let $G$ be an oriented graph with $n$ vertices and $\varepsilon n^2$ edges. We can assume that $\varepsilon < 1/4$ as otherwise the lemma follows easily by taking $c$ small enough and applying Corollary 3.2. Let $G' = (V, E)$ be the subgraph of $G$ obtained by removing every vertex with in-degree greater than $2\varepsilon n$. Observe that we removed at most $n/2$ vertices from $G$ and therefore $|V| \geq n/2$.

Let $U$ be a maximum acyclic set in $G'$, and assume that $|U| < \frac{\log n}{10\varepsilon \log (1/\varepsilon)}$. Since every vertex has in-degree at most $2\varepsilon n$, we have

$$\sum_{v \in U} d^-(v) < 2\varepsilon n |U|.$$

Let $R^*$ be the set of vertices in $V \setminus U$ with at least $5\varepsilon|U|$ out-neighbors in $U$, then $|R^*| \leq \frac{2n}{5}$, as otherwise the total number of edges directed into $U$ is more than $2\varepsilon n|U|$. Let $R = V \setminus (U \cup R^*)$ be the set of vertices outside $U$ with at most $5\varepsilon|U|$ out-neighbors in $U$, and we conclude that

$$|R| \geq n/2 - 2n/5 - |U| \geq n/20.$$

For every vertex $v \in R$ we define a set of vertices $S_v \subseteq U$ of size exactly $5\varepsilon|U|$ that contains $N^+(v) \cap U$. Using the inequality $\binom{n}{k} \leq (\frac{en}{k})^k$, we get that the total possible number of subsets of $U$ of this size is

$$\binom{|U|}{5\varepsilon|U|} \leq \left(\frac{e}{5\varepsilon}\right)^{5\varepsilon|U|} \leq 2^{\frac{\log n}{2\log(1/\varepsilon)} \cdot \log \frac{e}{5\varepsilon}} \leq n^{1/2}.$$

Therefore, by the pigeonhole principle, there is a set $R' \subset R$ of size at least $\frac{n}{20n^{1/2}} = \frac{n^{1/2}}{20}$ such that $|N^+(R') \cap U| \leq 5\varepsilon|U|$. Also, by Corollary 3.2, $R'$ contains an acyclic set $R'' \subseteq R'$ of size at least $\frac{1}{2} \cdot (\log n - 10)$. Note that $R'' \cup (U \setminus N^+(R'))$ is an acyclic set of size at least

$$|R''| + |U| - 5\varepsilon|U| \leq |U| + \frac{1}{2}(\log n - 10) - \frac{\log n}{2\log(1/\varepsilon)} \geq |U|,$$

assuming that $\varepsilon < 1/4$ and $n$ is large enough. Therefore, we get a contradiction as $U$ is not a maximum acyclic set, and the lemma follows. □

## 3.2 Acyclic colorings and coloring acyclic sets

In this subsection we give two building blocks that will be used later in the proof of Theorem 1. We first consider the case where we are given a coloring of the edges in some induced disjoint sets without a long monochromatic path, and we wish to color the edges between them while keeping the length of a longest monochromatic path relatively small. We then consider the case when we are given an acyclic set, and show how to color its edges with no long monochromatic path.

Given a set of vertices $A$ and a coloring $\varphi$ of $E(A)$, let $\ell_\varphi(A)$ be the length of a longest monochromatic path of in $A$ with respect to $\varphi$.

We prove the following.

**Lemma 3.4.** *Let $A_1, \ldots, A_k$ be disjoint sets, and let $\varphi$ be a coloring of $\bigcup_{i=1}^{k} E(A_i)$ with $(q+1)$ colors such that for every $1 \leq i \leq k$ we have $\ell_\varphi(A_i) \leq r$. Then $\varphi$ can be extended to a $(q+1)$ coloring $\varphi'$ of $E(\bigcup_{i=1}^{k} A_i)$ such that*

$$\ell_{\varphi'}\left(\bigcup_{i=1}^{k} A_i\right) \leq q(r+1) \cdot \lceil \sqrt[q]{k} \rceil.$$

*Proof.* Let $s$ be the minimal integer such that $k \leq s^q$. For each $1 \leq i \leq k$, we denote by $(i)$ the representation of $i$ in base-$s$ with exactly $q$ digits. That is, we represent each index as a vector of length $q$, letting $(i)_y$ be the $y$'th coordinate of $i$, and for each $1 \leq y \leq q$, $0 \leq (i)_y \leq s-1$.

We now define an acyclic coloring of the edges between the sets. For two sets $A_i$ and $A_j$, we color all the edges from $A_i$ to $A_j$ as follows. If there is an index $y$ such that $(i)_y < (j)_y$, we color all the edges in color $y$ (if more than one choice of $y$ is possible, choose one arbitrarily). Otherwise, we color all the edges from $A_i$ to $A_j$ in color $q+1$.

Note that by the definition of the new coloring, if a monochromatic path leaves a set $A_i$ it will never return to this set. Also, if $P$ is a monochromatic path colored by $1 \leq y \leq q$, then the $y$'th coordinates of all the sets visited by the path are all distinct, and therefore there are at most $s$ sets in the path. For a monochromatic path $P$ that is colored by $(q+1)$, if $A_j$ is visited by $P$ after it visits $A_i$ then $(i)_y \geq (j)_y$ for every coordinate $1 \leq y \leq q$, and there is a strict inequality in at least one coordinate. Therefore, the sums of the coordinates of all the sets visited by the path are all distinct, and hence such a path visits at most $sq$ distinct sets.

We conclude that every monochromatic path colored by $1 \leq y \leq q$ visits at most $s \leq sq$ distinct sets, and every monochromatic path colored by $q+1$ visits at most $sq$ sets.

A longest monochromatic path contains at most $r$ edges from each $A_i$, plus one edge that leaves $A_i$. It visits at most $qs = q\lceil \sqrt[q]{k} \rceil$ sets. Therefore, the length of this path is bounded by $q(r+1) \cdot \lceil \sqrt[q]{k} \rceil$, as claimed. □

In particular, we get the following immediate corollary, that generalizes an argument of Reimer [15].

**Corollary 3.5.** *Let $G$ be a $k$-colorable graph. There is a $(q+1)$-coloring of $E(G)$ with no monochromatic path longer than $q\lceil \sqrt[q]{k} \rceil$.*

We next prove similarly that one can color the edges of an acyclic graph $Z$ using $q$ colors with no monochromatic path longer than $\lceil \sqrt[q]{t} \rceil$, where $t$ is the length of a longest directed path in $Z$.

**Lemma 3.6.** *Let $Z$ be an acyclic graph in which a longest directed path has $t$ edges. There is a $q$-coloring of $E(Z)$ with no monochromatic path longer than $\lceil \sqrt[q]{t+1} \rceil$.*

*Proof.* Let $s$ be the minimal number such that $t + 1 \leq s^q$, and we prove that there a $q$-coloring with no monochromatic path longer than $s$. For a vertex $v$ let $\ell_Z(v)$ be the length of a longest path in $Z$ that ends in $v$, and for each $0 \leq i \leq t$ let

$$Z_i = \{v : \ell_Z(v) = i\}.$$

Note that these sets are well defined since $Z$ is an acyclic graph. Moreover, each $Z_i$ is an independent set and all the edges in the graph are directed from $Z_i$ to $Z_j$ for $j > i$. Let $(i)$ be the encoding of the number $i$ in base $s$ with exactly $q$ digits. Observe that for each $j > i$ there is an index $1 \leq y \leq q$ for which $(j)_y > (i)_y$. We now define the coloring as follows. For $j > i$, we color all the edges from $Z_i$ to $Z_j$ in color $y$ where $y$ is an index such that $(j)_y > (i)_y$. If there is more than one feasible choice of $y$, we choose one of them arbitrarily.

Again, we get that in every monochromatic path of color $y$, all the sets $Z_i$ that are visited by the path have distinct $y$'th coordinate, and there is at most one vertex from each $Z_i$. We conclude that every monochromatic path contains at most $s$ vertices, as required. □

We will also need the following well known and simple claim, and we give its proof for completeness.

**Claim 3.7.** *Let $G$ be an undirected graph with $m$ edges, than $G$ is $2\sqrt{m}$-colorable.*

*Proof.* Every optimal vertex coloring contains a vertex between every two color classes, hence $m \geq \binom{\chi(G)}{2}$ and we get that $\chi(G) \leq 2\sqrt{m}$ as required. □

**Proof of Theorem 1.** Let $G = (V, E)$ be a directed graph with $\frac{c_1 n^{2q}(\log n)^{1/q}}{(\log \log n)^{(q+2)/q}}$ edges. Define

$$X = \{v \in V : d^+(v) + d^-(v) \leq \left(\frac{n}{2q}\right)^q\}.$$

Note that $\Delta(G[X]) \leq (\frac{n}{2q})^q$ (when $G[X]$ is considered here as an undirected graph), and therefore, $\chi(G[X]) \leq (\frac{n}{2q})^q + 1$. By Corollary 3.5, there is a $(q+1)$-coloring of $E(X)$ with no monochromatic path longer than $n/2 + q$.

Let $Y = V \setminus X$, and $|Y| = m$. Since $|E(G)| \leq \frac{c_1 n^{2q}(\log n)^{1/q}}{(\log \log n)^{(q+2)/q}}$ and every vertex in $Y$ has at least $(\frac{n}{2q})^q$ incident edges, we have by double counting

$$m \leq \frac{2|E|}{(\frac{n}{2q})^q} \leq \frac{2c_1(2q)^q n^q (\log n)^{1/q}}{(\log \log n)^{(q+2)/q}}. \tag{3.1}$$

Consider the following procedure that partitions most of the vertices in $Y$ into $\Theta(\log \log n)$ families $Y^{(1)}, Y^{(2)}, \ldots$, where each family is composed of not too many acyclic sets.

Roughly speaking, we partition our vertices into acyclic sets as follows. We maintain an index $i$, and at the $i$'th step we find acyclic sets that cover approximately $m/2^i$ vertices, and group these sets to a *family* $Y^{(i)}$. To this end, at the beginning of each step we fix a number $a_i$ for which we are guaranteed by Lemma 3.3 that throughout the $i$'th step we can always find an acyclic set of size $a_i$. By the end of the procedure, most of the vertices are partitioned into families, and each family is composed of acyclic sets of the same size.

We initiate the procedure by taking $i, j = 1$ and $Y' = Y$, where $i$ represents the step number. At the beginning of step $i$, we fix $\varepsilon_i = \frac{|E(Y')|}{(m/2^i)^2}$ and $a_i = \frac{c \log(m/2^i)}{\varepsilon_i \log(1/\varepsilon_i)}$ (where $c$ is an absolute constant given by Lemma 3.3), and also set $j := 1$. The step ends when $|Y'| \leq m/2^i$.

At the $i$'th step, we repeatedly find an acyclic set $Z_{ij}$ in $Y'$ of size $a_i$, and set $Y' := Y' \setminus Z_{i,j}$ and then increase $j$ by one.

The procedure terminates when $|E(Y')| \leq \frac{n^{2q}}{(16q)^{2q}}$.

We first stress that for every $i, j$ we can find an acyclic set $Z_{ij}$ as claimed. Observe that during each step, the value $\varepsilon_i$ is fixed, the number of edges does not increase, and the number of vertices is at least $m/2^i$. Hence, $\varepsilon_i$ is an upper bound on edge density during the step. Therefore, Lemma 3.3 guarantees the existence of an acyclic set $Z_{ij}$ of size exactly $a_i$.

Let $Y^{(i)} = \{Z_{ij}\}$ be the family of acyclic sets that is constructed at step $i$, and let $k_i = |Y^{(i)}|$. We next bound the size of a longest monochromatic path in $\bigcup_{j=1}^{k_i} Z_{i,j}$ for every $i$, and start with $i = 1$.

Let $\varepsilon$ be the edge density in $Y$. Then $\varepsilon = \varepsilon_1/4$. It is easy to verify that $\frac{\log m - 1}{\log(1/\varepsilon) - 2} \geq \frac{\log m}{\log(1/\varepsilon)}$ (since clearly $m^2 \geq 1/\varepsilon$) and hence the size of each acyclic set is at least

$$\frac{c \log(m/2)}{\varepsilon_i \log(1/\varepsilon_i)} = \frac{c(\log m - 1)}{4\varepsilon(\log(1/\varepsilon) - 2)} \geq \frac{c \log m}{4\varepsilon \log(1/\varepsilon)}.$$

As we complete the first step just after we cover at least $m/2$ vertices with acyclic sets, we get that the number of acyclic sets satisfies $k_1 \leq \frac{(m/2) \cdot 4\varepsilon \log(1/\varepsilon)}{c \log m} + 1$. Moreover, since all acyclic sets are of the same size, the size of each $Z_{1j}$ is bounded by $\frac{m/2}{k_1 - 1}$.

By Lemma 3.6, we get that the edges of each acyclic set $E(Z_{1,j})$, for every $1 \leq j \leq k_1$, can be colored with $(q+1)$ colors with no monochromatic path longer than $\sqrt[q+1]{\frac{m/2}{k_1-1}} + 1$. We then can apply Lemma 3.4 to color the edges between $Z_{1,j_1}$ and $Z_{1,j_2}$ for each $j_1, j_2$. The number of sets is $k_1$, and the size of a longest monochromatic path in each of them is bounded by $\sqrt[q+1]{\frac{m/2}{k_1-1}} + 1$. We conclude that the edges of $\bigcup Z_{1,j}$ admit a $(q+1)$-coloring with no monochromatic path longer than

$$q(\sqrt[q]{k_1}+1) \cdot \left(\sqrt[q+1]{\frac{m/2}{k_1-1}} + 2\right) \leq 2q\left(\frac{m^q k_1}{2^q}\right)^{\frac{1}{q(q+1)}}.$$

Note that the edge density $\varepsilon$ in $Y$ satisfies

$$\varepsilon \leq \frac{|E|}{m^2} = \frac{c_1 n^{2q}(\log n)^{1/q}}{m^2(\log \log n)^{(q+2)/q}}. \tag{3.2}$$

Moreover, since we assume that the number edges is at least $\frac{n^{2q}}{(16q)^{2q}}$ (as otherwise the procedure terminates before this step begins), we get that $\frac{1}{\varepsilon} = \frac{m^2}{|E(Y)|} \leq \frac{m^2(16q)^{2q}}{n^{2q}}$ and therefore by (3.1)

$$\log(1/\varepsilon) \leq \log\left(\frac{(16q)^{2q}}{n^{2q}} \cdot \frac{4c_1^2(2q)^{2q}n^{2q}(\log n)^{2/q}}{(\log \log n)^{(2q+4)/q}}\right) \leq 2 \log \log n, \tag{3.3}$$

assuming that $n$ is large enough.

Therefore, we have the following bound.

$$\begin{aligned}
\frac{m^q k_1}{2^q} &\leq \frac{m^{q+1}\varepsilon \log(1/\varepsilon)}{2^{q-2}c \log m} \\
&\leq \frac{m^{q-1}c_1 n^{2q}(\log n)^{1/q} \log \log n}{2^{q-4}c \log m \cdot (\log \log n)^{(q+2)/q}} \\
&\leq \frac{8(2q)^{q(q-1)}c_1^q n^{q(q+1)} \log n \log \log n}{(\log \log n)^{q+2} c \log m} \\
&\leq \frac{8(2q)^{q(q-1)}c_1^q n^{q(q+1)}}{c(\log \log n)^{q+1}}
\end{aligned}$$

Where the first inequality follows from the bound on $k_1$, the second one from substituting $\varepsilon$ and $1/\varepsilon$ and using (3.2) and (3.3), and the third one one from substituting $m$ (and applying (3.1)). The last inequality follows from the bound $\log n \leq \log m$ that clearly holds.

Taking $c_1 < \frac{c^{1/q}}{8(2q)^q \cdot (16q^2)^{q+1}}$ (we do not try to optimize the dependence in $q$ here), we get that

$$2q\left(\frac{m^q k_1}{2^q}\right)^{\frac{1}{q(q+1)}} \leq \frac{1}{8q} \cdot \frac{n}{\sqrt[q]{\log \log n}}.$$

It is not difficult to verify that taking smaller values of $m$ only decreases the last expression. Hence, by repeating the same coloring method, we get that for each $i$, there is a $(q+1)$-coloring of the edges spanned by the union of the sets in $Y^{(i)}$ with no monochromatic path longer than $\frac{1}{8q} \cdot \frac{n}{\sqrt[q]{\log \log n}}$. Note that the number of vertices in $Y'$ is decreased by a factor of at least 2 in each step. Moreover, the procedure terminates when $|E(Y')| \leq \frac{n^{2q}}{(16q)^{2q}}$, and thus it essentially terminates if $|Y'| \leq \frac{n^q}{(16q)^q}$. Hence by applying (3.1) we get that the number of families is bounded by

$$\log \frac{(16q)^q m}{n^q} \leq \log \left( \frac{(16q)^q}{n^q} \cdot \frac{2c_1(2q)^q n^q (\log n)^{1/q}}{(\log \log n)^{(q+2)/q}} \right) \leq \log \log n,$$

assuming again that $n$ is large enough.

Let $\mathcal{W} = \bigcup_{i,j} Z_{ij}$ the set of vertices that were covered throughout the procedure. We color all the edges in $E(\mathcal{W})$ whose endpoints belong to sets from distinct $Y^{(i)}$'s with $(q+1)$ colors according to Lemma 3.4. Since the length of each monochromatic path spanned by the sets of each family is bounded by $\frac{1}{8q} \cdot \frac{n}{\sqrt[q]{\log \log n}}$, and the number of families is bounded by $\log \log n$, we get that the induced subgraph of all vertices that are contained in sets from $\bigcup Y^{(i)}$ admits a $(q+1)$-coloring with no path longer than

$$q \cdot \left( \frac{1}{8q} \cdot \frac{n}{\sqrt[q]{\log \log n}} + 1 \right) \cdot \left( \sqrt[q]{\log \log n} + 1 \right) \leq n/4.$$

Finally, after the procedure terminates, we are left with a set $Y'$ that satisfies $|E(Y')| \leq \frac{n^{2q}}{(16q)^{2q}}$. By Claim 3.7, we get that $G[Y']$ is $\frac{2n^q}{(16q)^q}$-colorable. Hence Corollary 3.5 implies that $E(Y')$ admits a $(q+1)$-coloring with no monochromatic path longer than $n/8 + 1$.

We therefore found a partition of $V(G)$ into three sets $X, Y', \mathcal{W}$. The set $E(X)$ admits a $(q+1)$-coloring with no monochromatic path longer than $n/2 + q$. The set $E(Y')$ admits a coloring with no monochromatic path longer than $n/8$, and $E(\mathcal{W})$ admits a coloring with no monochromatic path longer than $n/4$. We color all the edges that are either from $X$ to $Y'$ or from $X$ to $\mathcal{W}$ or from $Y'$ to $\mathcal{W}$ by the first color, and all the edges that are either from $Y'$ to $X$ or from $\mathcal{W}$ to $X$ or from $\mathcal{W}$ to $Y'$ by the second color. Every monochromatic path in $E(G)$ that leaves one of these three sets does not return to this set. Therefore, the length of a longest monochromatic path is bounded by the sum of the lengths of longest monochromatic paths in $X, Y', \mathcal{W}$, plus at most two edges between them. Thus $E(G)$ admits a coloring with no monochromatic path longer than $n/2 + q + n/8 + 1 + n/4 + 2 < n$, and Theorem 1 follows. □

**Proof of Proposition 4.** Let $G$ be a non-simple directed graph with $\left(\frac{n}{3q}\right)^{2q}$ edges. By Claim 3.7, $G$ is $2\left(\frac{n}{3q}\right)^q$-colorable, and hence by Corollary 3.5 it admits a $(q+1)$-coloring with no monochromatic path longer than $n$, and the proposition follows. □

## 4 Upper Bounds

In this section we provide an oriented graph for which every $q$-coloring of its edges contains a long monochromatic path. We start with the case of two colors. In Subsection 4.1 we define the notion of a $k$-pseudorandom digraph, and show that a random tournament on $n$ vertices is $\Theta(\log n)$-pseudorandom with high probability. We next show that every red-blue coloring of a $k$-pseudorandom digraph on $n$ vertices contains a directed red path of length $\Omega(\frac{n}{k})$ or a directed blue path of length $\Omega(n)$. This proves Theorem 3. We conclude this section by showing how to reduce the case of any fixed number of colors to the case of two colors, proving Theorem 2 and Proposition 5. For notational convenience we use $n$ in this section to denote the number of vertices in a Ramsey digraph $G$, rather than the length of a target path $P_n$ as in Theorem 2, Theorem 3 and Proposition 5.

## 4.1 Pseudorandom digraphs

We start with the following natural definition for $k$-pseudorandomness of directed graphs.

**Definition 4.1.** We say that a directed graph $G$ is *k-pseudorandom* if for every two disjoint sets $A, B$ such that $|A|, |B| \geq k$, there is at least one edge of $G$ from $A$ to $B$.

We first show the existence of a $\Theta(\log n)$-pseudorandom digraph by showing that a random tournament satisfies this property with high probability.

**Claim 4.2.** *A random tournament on $n$ vertices is $2 \log n$-pseudorandom with high probability.*

*Proof.* Let $T_n$ be a random tournament and fix two disjoint sets $A, B$ of size $k$ in $T_n$. Since every edge is oriented in each way uniformly and independently of the other edges, the probability that all the edges are directed from $B$ to $A$ is exactly $2^{-k^2}$. There are at most $\binom{n}{k}^2 \leq \frac{n^{2k}}{(k!)^2} \leq 2^{2k \cdot \log n - k \log k}$ choices of pairs of sets of size $k$, and thus by the union bound the probability that there are two sets of size $k$ for which all the edges are oriented in one direction is bounded by

$$2^{2k \cdot \log n - k^2 - k \log k} = o(1),$$

for $k = 2 \log n$. □

We also need the following property of $k$-pseudorandom digraphs.

**Claim 4.3.** *Let $G$ be a $k$-pseudorandom directed graph, and let $A_1, A_2, \ldots, A_t$ be disjoint sets, each of size at least $2k$. Then there is a directed path $v_1 v_2 \ldots v_t$, where for each $1 \leq i \leq t$, $v_i \in A_i$.*

*Proof.* We say that a vertex $u_j \in A_j$ is *good* if there is a directed path $u_j u_{j+1} \ldots u_t$ such that $u_s \in A_s$, $s = j, \ldots, t$. Clearly, our goal is to prove that there is a good vertex in $A_1$. Denote by $A_j^*$ the set of good vertices in $A_j$. By definition, every vertex in $A_t$ is good, and thus $|A_t^*| \geq 2k$. Also, if $u_{j+1} \in A_{j+1}^*$ and there is an edge from $u_j$ to $u_{j+1}$, then $u_j \in A_j^*$. Using a reverse induction, assume that we know that for some $j \leq t$, $|A_j^*| \geq k$, then since $G$ is $k$-pseudorandom all but at most $k$ of the vertices in $A_{j-1}$ have an edge to $A_j^*$. Therefore, all the vertices at most but $k$ in $A_{j-1}$ are actually in $A_{j-1}^*$, and thus $|A_{j-1}^*| \geq k$. We conclude that $|A_1^*| \geq k$ and thus there is a path as required. □

The following lemma shows that every $k$-pseudorandom graph contains a long path. The proof follows ideas from [7, 5].

**Lemma 4.4.** *Let $G$ a $k$-pseudorandom oriented graph on $n$ vertices. Then $G$ contains a directed path of length $n - 2k + 1$.*

*Proof.* Recall that the DFS (Depth First Search) is a graph algorithm that visits all the vertices of a (directed or undirected) graph $G$ as follows. It maintains three sets of vertices, letting $S$ be the set of vertices which we have completed exploring them, $T$ be the set of unvisited vertices, and $U = V(G) \setminus (S \cup T)$, where the vertices of $U$ are kept in a *stack* (a *last in, first out* data structure).

It is also assumed that some order $\sigma$ on the vertices of $G$ is fixed, and the algorithm prioritizes vertices according to $\sigma$. The DFS starts with $S = T = \emptyset$ and $U = \{\sigma_1\}$.

While there is a vertex in $V(G) \setminus S$, if $U$ is non-empty, let $v$ be the last vertex that was added to $U$. If $v$ has a neighbor $u \in T$, the algorithm inserts $u$ to $U$ and repeats this step. If $v$ does not have a neighbor in $T$ then $v$ is popped out from $U$ and is inserted to $S$. If $U$ is empty, the algorithm chooses an arbitrary vertex from $T$ and pushes it to $U$.

We are now proceed to the proof of the lemma. We Execute the DFS algorithm. We let again $S, T, U$ be three sets of vertices as defined above. At the beginning of the algorithm, all the vertices are in $T$, and at each step a single vertex either moves from $T$ to $U$ or from $U$ to $S$. At the end of the algorithm, all the vertices are in $S$. Therefore, at some point we have $|S| = |T|$. Observe crucially that all the vertices in $U$ form a directed path, and that there are no edges from $S$ to $T$. We conclude that $|S| \leq k - 1$, and therefore $|U| \geq n - 2k + 2$, so there is a directed path with $n - 2k + 1$ edges in $U$, as required. □

## 4.2 The case of two colors

In this subsection we prove that every $k$-pseudorandom directed graph on $n$ vertices has a monochromatic red path of length $\Omega(\frac{n}{k})$ or a monochromatic blue path of length $\Omega(n)$, and this will prove Theorem 3. We prove the following main lemma.

**Lemma 4.5.** *Let $G$ be a $k$-pseudorandom directed graph on $n$ vertices. Then every red-blue coloring of its edges yields a red directed path of length $\frac{n}{28k}$ or a blue directed path of length $n/28$.*

*Proof.* Fix a red-blue coloring of $E(G)$. Let $G_R$ be the *red graph*, that contains only the red edges.

If $G_R$ contains a directed path of length $\frac{n}{14k}$, we are done. Otherwise, by the Gallai-Roy theorem (Theorem 2.2), $G_R$ is $\frac{n}{14k}$-colorable and therefore has a partition into $\frac{n}{14k}$ independent sets. We partition these independent sets into sets of size exactly $7k$, rounding down the remaining vertices (in particular, we remove every independent set smaller than $7k$). Since we remove at most $7k \cdot \frac{n}{14k} = n/2$ vertices, we remain with $t \geq \frac{n}{14k}$ independent sets $B_1, B_2, \ldots B_t$, each containing exactly $7k$ vertices.

Note that each $B_i$, $1 \leq i \leq t$, spans a $k$-pseudorandom graph and contains only blue edges. Therefore, by Lemma 4.4, each $B_i$ contains a blue path of length at least $5k$. Since there is a directed edge from the last $k$ vertices of this path to the first $k$ vertices in the path, we conclude that each $B_i$ contains a directed blue cycle $C_i$ of length at least $3k$.

We next construct an auxiliary graph $H$ on $t$ vertices, each vertex corresponds to a cycle $C_i$, and with a slight abuse of notation we denote these vertices by $C_1, C_2, \ldots, C_t$. The graph $H$ is a complete symmetric directed graph, and we color the edges from $C_i$ to $C_j$ by blue if there are at least $k$ vertices in $C_i$ that have blue edges directed towards $C_j$, and red otherwise. Since $H$ is complete and symmetric, by Raynaud's theorem (Theorem 2.1), it has a monochromatic path of length $t/2 \geq \frac{n}{28k}$. Let $C_{i_1}, C_{i_2}, \ldots, C_{i_{t/2}}$ be the vertices in $H$ along the path, each such vertex represents a cycle. We complete the proof by considering the following two cases.

**$H$ contains a red path of length $t/2$.** For each $C_{i_j}$ there is a set $R_{i_j}$ for which only red edges are going towards $C_{i_{j+1}}$. Moreover, for every $1 \leq j \leq t/2$ we have $|R_{i_j}| \geq 2k$. Observe that all the edges from $R_{i_j}$ to $R_{i_{j+1}}$ are red, and therefore by Claim 4.3, there is a red path of length $t/2 \geq \frac{n}{28k}$ with exactly one vertex from each $R_{i_j}$, as required.

**$H$ contains a blue path of length $t/2$.** Call a vertex in $C_{i_j}$ an *endpoint* if it has a blue edge towards $C_{i_{j+1}}$. By the assumption, there are at least $k$ endpoints in $C_{i_j}$ for every $1 \le j \le t/2$, and therefore each vertex in $C_{i_j}$ there is a path of length at least $k-1$ that ends at some endpoint, in which we can travel along one additional edge towards $C_{i_{j+1}}$. We construct a blue path of length $\frac{n}{28}$ by taking an arbitrary path of length $k-1$ that ends at some endpoint in $C_{i_1}$, moving through the endpoint to $C_{i_2}$, walking through such a path to an endpoint of $C_{i_2}$ and so on till we arrive at $C_{i_{t/2}}$, where we can again walk along a path of length at least $3k-1$ (that visits all the vertices in $C_{i_{t/2}}$). We conclude that there is a blue path of length at least $k \times \frac{n}{28k} = \frac{n}{28}$, as claimed.

Lemma 4.5 and Theorem 3 follow. □

**Explicit constructions.** Given an explicit construction of a $k$-pseudorandom tournament on $n$ vertices, our approach shows that every red-blue coloring of such a tournament has a monochromatic path of length $\Omega(\frac{n}{k})$. To the best of our knowledge, the best construction of a $k$-pseudorandom tournament is given by Quadratic Residue tournaments, defined as follows (see [4], Chapter 9). Let $p \equiv 3 \mod 4$ be a prime. The vertices of the tournament $T_p$ are all the elements in the finite field $\mathbb{Z}_p$. For two vertices $i$ and $j$, there is an edge from $i$ to $j$ if and only if $i-j$ is a quadratic residue. It can be shown that since $p \equiv 3 \mod 4$, $-1$ is a quadratic nonresidue and therefore for each $i,j$ there either there is an edge from $i$ to $j$ or an edge from $j$ to $i$ but not both. This construction gives $k = \Theta(\sqrt{n})$ (see [2] and [4]).

However, every red-blue coloring of any tournament yields a monochromatic path of length $\sqrt{n} - 1$. To prove this, observe that if both the red and blue graphs, formed respectively by taking all red and blue edges in some fixed coloring do not have a path of length $\sqrt{n} - 1$, then by the Gallai-Roy theorem (Theorem 2.2) they are both $(\sqrt{n}-1)$-colorable. It is easy and well known (see, e.g., [13], Chapter 9, Problem 3) that given two graphs $G_1, G_2$ we have $\chi(G_1 \cup G_2) \le \chi(G_1) \cdot \chi(G_2)$. We conclude that $K_n$, which is the union of the red and blue graph is $(n-1)$-colorable, and we get a contradiction.

Hence, only explicit constructions of $o(\sqrt{n})$-pseudorandom graphs will be interesting for our problem.

## 4.3 The general case

Here we prove by induction on $q$ that for every $k$-pseudorandom directed graph $G$ on $28kn^q$ vertices we have $G \to (\overrightarrow{P_n})_{q+1}$. Theorem 2 will clearly follow by the fact that a random tournament is $\Theta(\log n)$-pseudorandom with high probability (Claim 4.2).

The base case of the induction ($q = 1$) follows directly from Lemma 4.5. Suppose that the result holds for $q$ colors and we next prove it for $q+1$ colors. Indeed, let $G$ be a $k$-pseudorandom graph on $28kn^q$ vertices. Fix a coloring of $E(G)$ with the colors $1, 2, \ldots, q+1$. Denote by $G_{q+1} \subseteq G$ the subgraph with all edges that are colored $q+1$. If $G_{q+1}$ contains a monochromatic path of length $n$ we are done. Otherwise, by the Gallai-Roy Theorem (Theorem 2.2), we get that $G_{q+1}$ is $n$-colorable and hence contains an independent set $A$ of size $\frac{28kn^q}{n} = 28kn^{q-1}$.

Note that $A$ spans a subgraph of $G$ and hence $G[A]$ is $k$-pseudorandom. Also, the colors of all the edges of $G[A]$ are among $\{1, 2, \ldots, q\}$. Therefore, by the induction hypothesis $G[A] \to (\overrightarrow{P_n})_q$, and we conclude that $G$ contains a monochromatic path of length $n$, as desired. □

**Proof of Proposition 5.** The proof for non-simple directed graphs follows similar lines. Note that the Gallai-Roy theorem is valid for non-simple directed graphs as well. Therefore, we can use the same induction on $q$. The base case ($q = 1$) follows from Raynaud's theorem (Theorem 2.1). Suppose that the result holds for $q$ colors, the correctness for $q + 1$ colors follows by taking a complete symmetric graph on $n^q$ vertices, and considering the subgraph with all edges colored by the $(q+1)$'th color. If this graph contains a directed path of length $n$ we are done, otherwise we find an induced subgraph of order $n^{q-1}$ in which all edges are colored $1, 2, \ldots, q$, and applying the induction hypothesis, Proposition 5 follows. □

## 5 Concluding remarks and open problems

We proved nearly tight bounds for the size Ramsey number of a directed path for oriented graphs. We proved that every red-blue coloring of the edges of a $k$-pseudorandom graph on $n$ vertices contains a red path of length $\Omega(\frac{n}{k})$ or a blue path of length $\Omega(n)$, but it might be the case that this approach can also give better symmetric Ramsey bounds. An interesting question is whether every red-blue coloring of a $k$-pseudorandom graph contains a monochromatic path of length $\Omega(\frac{n}{\sqrt{k}})$. Clearly every progress in this direction will improve our upper bounds.

Another related question is about the asymptotic behavior of the maximum length of a monochromatic path in every red-blue coloring of a random tournament. Here we proved that every tournament $T$ has a coloring with no monochromatic path longer than $O(\frac{n}{\sqrt{\log n}})$, and also that with high probability a random tournament $T_n$ has a monochromatic path of length $\Omega(\frac{n}{\log n})$ in every red-blue coloring. It would be very interesting to close the gap between these bounds.

When proving the lower bound on the size Ramsey number, we study the minimal number $k$ for which a certain graph can be partitioned into $k$ acyclic sets. This parameter was studied, e.g., in [1], and it was conjectured that in every oriented graph $G = (V, E)$ there is an acyclic set of size $(1 + o(1))\frac{|V|^2}{|E|} \cdot \log \frac{|E|}{|V|}$. If this conjecture is true, then our lower bound can be slightly improved.

It is easy to verify that there is no $k$-pseudorandom oriented graph on $n$ vertices for $k \leq \frac{\log n}{2}$, as every such graph has an acyclic set of size $\log n$ and therefore has two sets of size $\frac{\log n}{2}$ with no edges in one of the directions. On the other hand we proved that for $k = 2 \log n$ such graphs do exist. Hence, it will be interesting to determine the minimum $k$ for which there is a $k$-pseudorandom oriented graph on $n$ vertices. Another appealing question is to provide better explicit constructions of $k$-pseudorandom oriented graphs.

## References

[1] R. Aharoni, E. Berger and O. Kfir, Acyclic systems of representatives and acyclic colorings of digraphs, *J. of Graph Theory* 59 (2008), 177–189.

[2] N. Alon, Eigenvalues, geometric expanders, sorting in rounds and Ramsey theory, *Combinatorica* 6 (1986), 207–219.

[3] N. Alon and F. R. K. Chung, Explicit construction of linear sized tolerant networks, *Discrete Math.* 72 (1988), 15–19.

[4] N. Alon and J. Spencer, **The probabilistic method**, Third edition, Wiley, 2008.

[5] J. Bang-Jensen and S. Brandt, Expansion and Hamiltonicity in Digraphs, Manuscript.

[6] J. Beck, On size Ramsey number of paths, trees and circuits. I, *J. Graph Theory* 7 (1983), 115–129.

[7] S. Brandt, H. Broersma, R. Diestel and M. Kriesell, Global connectivity and expansion: long cycles in $f$-connected graphs, *Combinatorica* 26 (2006), 17–36.

[8] P. Erdős, R. J. Faudree, C. C. Rousseau and R.C. Schelp, The size Ramsey number, *Period. Math. Hungar.* 9 (1978), 145–161.

[9] P. Erdős and E. Szemerédi, On a Ramsey type theorem, Collection of articles dedicated to the memory of Alfréd Rényi, I., *Period. Math. Hungar.* 2 (1972), 295–299.

[10] J. Friedman and N. Pippenger, Expanding graphs contain all small trees, *Combinatorica* 7 (1987), 71–76.

[11] T. Gallai, On directed paths and circuits, in **Theory of graphs**, edited by P. Erdős and G. Katona, Academic Press, 1968, 115–118.

[12] A. Gyárfás, Vertex coverings by monochromatic paths and cycles, *J. Graph Theory* 7 (1983), 131–135.

[13] L. Lovász, **Combinatorial problems and exercises**, Second edition, AMS publishing, 2007.

[14] H. Raynaud, Sur le circuit hamiltonien bi-coloré dans les graphes orientés, *Period. Math. Hungar.* 3 (1973), 289–297.

[15] D. Reimer, The Ramsey size number of dipaths, *Discrete Math.* 257 (2002), 173–175.

[16] B. Roy, Nombre chromatique et plus longs chemins d'un graphs, *Rev. AFIRO* 1 (1967), 129–132.